\documentclass{article}

\usepackage{amsmath}
\usepackage{amssymb}
\usepackage{amsthm}

\theoremstyle{plain}
\newtheorem{theorem}{Theorem}[section]
\newtheorem{lemma}[theorem]{Lemma}
\newtheorem{definition}[theorem]{Definition}

\newtheorem{proposition}[theorem]{Proposition}

\def \vect#1 {\overrightarrow{#1}}

\def \Cc {{\cal C}}

\def \Qq {{\cal Q}}

\def \kbar {\bar k}

\def \Gal {{\rm Gal}\,}

\def \BBB {\mathbf{B}}
\def \OOO {\mathbf{O}}
\def \SO {\mathbf{SO}}

\def \Div {{\rm Div}\,}
\def \divi {{\rm div}\,}

\def \Pic {{\rm Pic}\,}

\def \rank {{\rm rank}\,}

\def \Divg {\Div ' \Cc}
\def \Picg {\Pic ' \Cc}
\def \fim  {\varphi_\lambda}
\def \Om  {\OOO _3 (\Omega)}
\def \SOm  {\SO _3 (\Omega)}

\def\harr#1#2{\smash{\mathop{\hbox to .5in{\rightarrowfill}}
\limits^{\scriptstyle#1}_{\scriptstyle#2}}}

\def\varr#1#2{\llap{$\scriptstyle #1$}\left\downarrow \vcenter to
.5in{}\right.\rlap{$\scriptstyle #2$}}

\begin{document}
\title{About the equivalence of divisor classes
on hyperelliptic curves and a quotient of linear forms by an
algebraic group action}
\author{V.G. Lopez Neumann}

\date{}

\maketitle

\section{introduction}
For a hyperelliptic curve of genus $g$, a divisor in general
position of degree $g+1$ is given by polynomial equations. There
is an action from an algebraic group on the representations of
divisors by polynomials which fixes divisor classes. This
structure reduces the question of rationality of divisor classes
to rationality of polynomials which is more easy to control. This
structure was inspired by \cite{CF}.


\section{notation}
Let $k$ be a perfect field of characteristic different from $2$,
$\kbar$ the algebraic closure of $k$ and $G=\Gal (\kbar / k)$ its
Galois group. If $A$ is a set where $G$ acts, we call $A(k)$ the set
of elements invariants by Galois group, they are the rational points
or rational elements of $A$. A map is called rational if it commutes
with the Galois group.

Let~$\Cc$ be a hyperelliptic curve of genus~$g$ defined over~$k$
by the affine equation
$$
Y^2 = F(X) = f_0 + f_1 X + \cdots + f_{2g+2}X^{2g+2} \in k[X],
$$
where $F$ has no multiple roots. The points with $1/X=0$,
$Y/X^{g+1}=\pm\sqrt{f_{2g+2}}$ on the completion are called the
points at infinity of $\Cc$. We denote them $\infty^\pm$.

\begin{definition}
An effective divisor~$D$ on~$\Cc$ of degree~$g+1$ is in general
position if it is given by
$$
U(X)=0, \quad Y=W(X),
$$
where~$U, \, W \in \kbar[X]$ have degree at most~$g+1$, such that
there exists a polynomial $V \in \kbar[X]$ for which
$$
F=W^2-UV .
$$
If~$U$ has degree less than~$g+1$, then some of the points of the
support of~$D$ are at infinity.
\end{definition}
For a divisor $D$, we write $\overline{D}$ its conjugate by the
involution which sends a point $(X,Y)$ into the point $(X,-Y)$, this
will be called the $\pm Y$ involution.


For simplicity, let $\Divg$ be the set of effective divisors~$D$
on~$\Cc$ of degree~$g+1$ in general position, and $\Picg$ be the
set of linear equivalence classes of such divisors. An element of
$\Picg (k)$ not necessarily comes from $\Divg (k)$ (see
\cite{CM}).

Let $(\kbar ^{g+2})^*$ be the $\kbar$-vector space of linear forms
on $g+2$ variables, denoted by $x=(x_0, x_1, \ldots , x_{g+1})$.
Let $\lambda:(\kbar ^{g+2})^* \longrightarrow \kbar[X]_{\le g+2}$
be the linear map given by
$$
x_i \mathop{\longrightarrow}^\lambda X^i \quad \text{for} \quad 0
\le 1 \le g+1 .
$$
Let $L (\Cc)$ be the subset of $((\kbar ^{g+2})^*)^3$ given by
$$
t= \left(
   \begin{array}{c}
   u \\ v \\ w
   \end{array}
   \right)
 \in L(\Cc) \Longleftrightarrow
F=\lambda(w)^2 -\lambda(u)\lambda(v)
$$
Then $\lambda$ extends to a map $\fim$ from $L(\Cc)$ to $\Divg$,
where $t \in L(\Cc)$ is sent to the divisor given by
$$
\lambda(u)(X)=0 , \quad Y=\lambda(w)(X) .
$$

Our goal is to find an algebraic group $\OOO$ acting on $L(\Cc)$
such that $\lambda$ induces a rational bijection between $L(\Cc) /
\OOO$ and $\Picg$.

Let $\Om$ be the orthogonal group relative to $\Omega$, which is
given by:
$$
\Om =
 \left\{
 A \in \text{Mat} \, (3,\kbar) \mid A^* \Omega A = \Omega
 \right\} ,
$$
where $A^*$ is the transpose of $A$ and $\Omega$ is the matrix
$$
\Omega = \left(
\begin{array}{rrr}
0    & -1/2 & 0 \\
-1/2 & 0    & 0 \\
0    & 0    & 1
\end{array}
\right) .
$$
If $A \in \Om$ then $\det A = \pm 1$. Call $\SOm$ the subgroup of
$\Om$ given by the additional condition $\det A = 1$.


\section{The action of $\SOm$ on linear forms}
\begin{definition}
Let $t \in L(\Cc)$. We have an action from the group $\Om$ on
$L(\Cc)$ given by matrix multiplication
$$
A \cdot t = A t \quad \text{where } A \in \Om .
$$
\end{definition}


Let $\BBB$ be the subgroup of $\SOm$ generated by the elements
$$
\left(
\begin{array}{rrr}
a & 0   & 0 \\
0 & 1/a & 0   \\
0 & 0   & 1
\end{array}
\right) \, , \, \left(
\begin{array}{rrr}
1    & 0 & 0   \\
b^2 & 1 & -2b \\
-b    & 0 & 1
\end{array}
\right) ,
$$
where $a \in \kbar^*$ and $b \in \kbar$.

\begin{lemma}\label{bijdiv}
The map $\fim : L(\Cc) \longrightarrow \Divg$ induces a rational
bijection between $L(\Cc) / \BBB$ and $\Divg$.
\end{lemma}
\begin{proof}
Let $D \in \Divg$ be given by
$$
U(X)=0 \, , \quad Y=W(X) .
$$
Any other representation of $D$ will be given by
$$
\begin{array}{rl}
a \, U(X)= 0 \, , & Y=W(X) , \\
     U(X)=0  \, , & Y=W(X) + b \, U(X) ,
\end{array}
$$
or a combination of them, where $a \in \kbar^*$ and $b \in \kbar$.
The conclusion is clear.
\end{proof}

By definition of divisors in general position, the map $[\fim] :
L(\Cc) \longrightarrow \Picg$ is surjective. The next proposition
proofs that the linear equivalence classes of effective divisors
of degree~$g+1$ in general position are in~$1$ to~$1$
correspondence with proper equivalence classes of representations
of $F$ by $W^2-UV$, where $U,\, V,\, W$ have degree at most $g+1$.
The involution $Y \mapsto -Y$ corresponds to $U\mapsto U$, $V
\mapsto V$, $W \mapsto -W$. Hence an improper automorphism of
$W^2-UV$ (i.e. of determinant $-1$) takes the corresponding
element of $\Picg$ into its conjugate.
\begin{proposition}\label{Cassels}
The map $L(\Cc) / \SOm  \longrightarrow \Picg$ induced by $[\fim]$
is a rational bijection. An element from $\Om$ not in $\SOm$ sends
a divisor class into its conjugate by $\pm Y$ involution, i.e. for
an element $A \in \Om$ such that $\det A = -1$, we have $[\fim]
(At) = [\overline{D}]$, where $t \in L(\Cc)$ and $D= \fim(t) \in
\Divg$.
\end{proposition}
\begin{proof}
Suppose we are proved the first part of the proposition. The last
part is then clear: let $A \in \Om$ such that $\det A = -1$ and
let $\varepsilon \in \Om$ be the element sending $(u,v,w)$ into
$(u,v,-w)$; in particular $A \varepsilon \in \SOm$. For $t \in
L(\Cc)$ we have:
$$
[\fim](At)= [\fim] (A \varepsilon (\varepsilon t)) = [\fim]
(\varepsilon t)= \overline{[\fim](t)} .
$$

We only need to prove that for all $t_1,t_2 \in L (\Cc)$, we have:
$$
[\fim](t_1)=[\fim](t_2) \Longleftrightarrow \exists A \in \SOm,
\text{ such that } t_2= A \, t_1 .
$$

For $i=1,2$, let $D_i=\fim(t_i)$ and $U_i=\lambda(u_i),\,
V_i=\lambda(v_i),\, W_i=\lambda(w_i)$ the corresponding
polynomials, where
$$
t_i=
 \left(
 \begin{matrix}
 u_i \\ v_i \\ w_i
 \end{matrix}
 \right) .
$$

Suppose that~$D_1 \sim D_2$. Let $D_1'$ be the divisor such that
$$
\divi \left( \frac{Y-W_1(X)}{V_1(X)} \right) = D_1 - D_1' ,
$$
then $D_1' \sim D_2$. By the Riemann-Roch theorem, the functions
$$
1 \, , \,   \frac{Y-W_1(X)}{V_1(X)}
$$
form a basis of the vector space $\{ f \in \kbar(\Cc) \mid D_1' +
\divi f \succ 0 \} \cup \{ 0 \}$. If $D_1' \neq D_2$ there exists $a
\in \kbar$ such that
$$
\divi \left( \frac{Y-W_1(X)}{V_1(X)} + a \right) = D_2 - D_1' .
$$
The divisor~$D_2$ is then given by
\begin{equation}\label{represent}
U_1(X) + a^2 V_1(X) - 2aW_1(X)=0 , \quad  Y=W_1(X) - a V_1(X) ,
\end{equation}
By lemma \ref{bijdiv}, there exists $B \in \BBB$ such that $D_2=\fim
(B t_2)$ is exactly represented by (\ref{represent}). Then, we may
suppose
$$
U_2 = U_1(X) + a^2 V_1(X) - 2aW_1(X) \quad \text{and} \quad W_2 =
W_1(X) - a V_1(X) .
$$
The relation on linear forms is $t_2 = A \, t_1$, where
$$
A=\left(
\begin{array}{rrr}
1 & a^2 & -2a \\
0 & 1   & 0   \\
0 & -a  & 1
\end{array}
\right) \in \SOm.
$$
If $D_1' = D_2$, then $D_2$ is given by
$$
D_1'=D_2 :\, V_1(X)=0, \, Y=-W_1(X) .
$$
In this case, by lemma \ref{bijdiv}, we may suppose $U_2=V_1$ and
$W_2=-W_1$. The relation $t_2 = A \, t_1$ is then given by
$$
A=\left(
\begin{array}{rrr}
0 & 1 & 0 \\
1 & 0 & 0   \\
0 & 0 & -1
\end{array}
\right) \in \SOm.
$$
In the other sense, suppose that there exists $A \in \SOm)$ such
that $t_2 = A t_1$.


Working explicitly with the expression $A^* \, \Omega \, A = \Omega$
and $\det A=1$, we find that $\SOm$ is generated by
$$
\left(
\begin{array}{rrr}
a & 0   & 0 \\
0 & 1/a & 0   \\
0 & 0   & 1
\end{array}
\right) \, , \, \left(
\begin{array}{rrr}
0 & 1 & 0   \\
1 & b^2 & 2b \\
0 & -b & -1
\end{array}
\right) ,
$$
where $a \in \kbar^*$ and $b \in \kbar$. For each of them we have
$\fim (t) = \fim (A t)$, for all $t \in L(\Cc)$.
\end{proof}


\section{Quadratic forms associated to divisor classes}

In this section, we will obtain an easy method to know when a given
divisor in general position of degree $g+1$ doesn't represents a
rational divisor class.

Let $\Qq(\kbar^{g+2})$ be the set of quadratic forms on $g+2$
variables. The map $\lambda$ induces a map from $\Qq(\kbar^{g+2})$
to $\kbar[X]$. We are concerned in the inverse image of $F$. Let
$\lambda^{-1}(F)$ be the subset of elements $S \in
\Qq(\kbar^{g+2})$ such that $\lambda(S)=F$.

We obtain a map $q_{\Omega} : L(\Cc) \longrightarrow
\Qq(\kbar^{g+2})$ given by:
$$
q_{\Omega} (t) = t^* \, \Omega \, t = w^2 -uv .
$$
A quadratic form on the image of $q_{\Omega}$ is clearly of rang
$2$ or $3$ and it lives on $\lambda^{-1}(F)$. Let $\Qq(\Cc)$ be
the subset of $\lambda^{-1}(F)$ of qudratic forms of rang less or
equal than $3$. In fact $\Qq(\Cc)$ is the image of $L(\Cc)$ by
$q_{\Omega}$.

\begin{lemma}\label{bilinear}
The map from $L(\Cc) / \Om$ to $\Qq(\Cc)$ is a rational bijection
induced by $q_{\Omega}$. By abuse of notation, we call this map by
$q_{\Omega}$ too.
\end{lemma}
\begin{proof}
It is clear that $q_{\Omega} : L(\Cc) / \Om \longrightarrow
\Qq(\Cc)$ is a well defined surjective rational map. We only need
to prove the injectivity.

Pick $t_1,t_2 \in L(\Cc)$ such that
$S=q_{\Omega}(t_1)=q_{\Omega}(t_2) \in \Qq(\Cc)$. For $i=1,2$, we
write
$$
t_i= \left(
     \begin{array}{c}
     u_i \\ v_i \\ w_i
     \end{array}
     \right) .
$$
Let $Q$ be the bilinear form associated to $S$. We have by
definition
$$
\text{rad} \, (S)= \{ \vect{x} \in \bar{k}^{g+2} \mid Q(\vect{x} ,
\vect{y} )=0,
             \, \forall \vect{y} \in \bar{k}^{g+2} \},
$$
where
$$
Q( \vect{x} , \vect{y} ) = w_i( \vect{x} ) w_i( \vect{y} ) -
\frac{1}{2} \left( u_i( \vect{x} ) v_i( \vect{y} ) + u_i( \vect{y} )
v_i( \vect{x} ) \right) , \quad \text{for } i=1,2 .
$$
It is clear that $Z(u_i,v_i,w_i)$, the set of zeros of the linear
forms~$u_i,\, v_i,\, w_i$ is included in~$\text{rad} \, (S)$. On the
other hand, we have $\rank S =3$ if and only if $\dim \langle
u_i,v_i,w_i \rangle =3$
 by linear algebra. Finally $\rank S >1$, because if not:
$$
z^2=w_i^2-u_i v_i \quad \text{and} \quad F=\lambda(z)^2;
$$
but $F$ has no multiple factors. That means $\rank S = \dim
\langle u_i,v_i,w_i \rangle$.

We conclude that $\text{rad} \, (S)=Z(u_i,v_i,w_i)$ for $i=1,2$.
Then
\begin{equation}\label{generator}
q_{\Omega}(t_1)=q_{\Omega}(t_2) \Longleftrightarrow \langle
u_1,v_1,w_1 \rangle = \langle u_2,v_2,w_2 \rangle .
\end{equation}
If $\rank S=3$, there exists a matrix $A$ such that
$$
t_2 = A t_1 \quad \text{and} \quad t_1^* A^* \Omega A t_1 = t_1^*
\Omega t_1 .
$$
By linear independence of $u_1,v_1,w_1$, we have $A \in \Om$.

If $\rank S=2$, there exists independents linear forms $u,v$ such
that $S=-uv$. We can suppose then that $u_2=u$, $v_2=v$ and $w_2=0$.
By (\ref{generator}), there exists $a,b,c,d,e,f \in \kbar$ such that
for any $x,y,z \in \kbar$, we have
$$
\left(
\begin{array}{r}
u_1 \\
v_1   \\
w_1
\end{array}
\right) = \underbrace{\left(
\begin{array}{rrr}
a & b & x \\
c & d & y   \\
e & f & z
\end{array}
\right)}_A \cdot \underbrace{\left(
\begin{array}{r}
u \\
v   \\
0
\end{array}
\right)}_t .
$$
Independence of $u,v$ and the relation $w_1^2 - u_1 v_1 = -uv$ give:
\begin{equation}\label{constantrel}
\begin{aligned}
e^2&=ac , \\
f^2&=bd , \\
2ef&=ad+bc-1.
\end{aligned}
\end{equation}
Replacing these relations on $A^* \Omega A = \Omega$, we obtain the
system on $x,y,z$:
$$
\begin{aligned}
cx+ay&=2ez , \\
dx+by&=2fz ,\\
xy+1 &=z^2 .
\end{aligned}
$$
On picking
$$
\begin{aligned}
x&=2(af-eb) , \\
y&=2(de-cf) ,\\
z&=ad-bc ,
\end{aligned}
$$
and taking into account the relations on (\ref{constantrel}), we
obtain a matrix $A \in \Om$.
\end{proof}

By proposition \ref{Cassels}, we have the commutative rational maps:
$$
\begin{array}{ccc}
L(\Cc)/\SOm & \harr{[\fim]}{\cong} & \Picg \\
\varr{}{} & & \varr{}{} \\
L(\Cc)/\Om & \harr{[\fim]}{\cong} & \Picg / \pm Y \\
\varr{q_{\Omega}}{\cong} & & \\
\Qq(\Cc) & &
\end{array}
$$
which gives us a double cover from $\Picg$ on $\Qq(\Cc)$.

In particular, this means that a divisor class $[D]$ modulo $\pm Y$
involution is rational if and only if the corresponding quadratic
form $S=q_{\Omega} \circ [\fim]^{-1} ([D])$ is rational. Attention,
if $S$ is rational, it doesn't implies that the divisor class $[D]$
is rational.

\vskip 1cm

\thanks{
\noindent V.G. Lopez Neumann \\
        UFMG, Belo Horizonte, MG, Brazil,\\
        supported by CNPq - Brazil\\
        email gonzalo@mat.ufmg.br}


\begin{thebibliography}{CF}
\bibitem[CF]{CF} {\sc J.\,W.\,S. Cassels \& E.\,V. Flynn},
{\it Prolegomena to a Middlebrow Arithmetic of Curves of Genus~$2$}.
London Math. Soc. Lecture Note Series 230, Cambridge, 1996.

\bibitem[CM]{CM} Coray, D., Manoil, C.:
On large Picard groups and the Hasse Principle for curves and K3
surfaces. Acta Arith. {\bf 76}, 165--189 (1996)


\end{thebibliography}
\end{document}